\documentclass{amsart}
\usepackage{amscd}
\theoremstyle{plain}
\newtheorem{Thm}{Theorem}[section]
\newtheorem{Prop}[Thm]{Proposition}
\newtheorem{Cor}[Thm]{Corollary}
\newtheorem{Lemm}[Thm]{Lemma}
\newtheorem{ThmA}{Theorem A}

\newtheorem{ThmB}{Theorem B}

\newtheorem{ThmC}{Theorem C}

\newtheorem{ThmD}{Theorem D}

\newtheorem{prfB}{\bf Proof$\!\!$}
\theoremstyle{remark}
\def\eps{\epsilon}
\newtheorem{Defn%
}[Thm]{\bf Definition}
\newtheorem{Term%
}[Thm]{\bf Terminology}
\newtheorem{Rmk%
}{\bf Remark}

\newtheorem{Notn%
}[Thm]{\bf Notation}

\newcommand{\abs}[1]{\lvert#1\rvert}
\newcommand{\norm}[1]{\lVert#1\rVert}
\newcommand{\hardhyphen }%
{\nobreakdash-\hspace{0pt}}
\begin{document}

\author[C. Pugh]{Charles Pugh}
\address{Mathematics Department
  \\ University of California \\
Berkeley California, 94720, U.S.A.}
  \email{pugh@math.berkeley.edu}
\author[M. Shub]{Michael Shub}
\address{Watson Research
Center
  \\ IBM \\ Yorktown Heights NY,
10598, USA}
  \email{shub@us.ibm.com}
 \thanks{Shub's work was partly funded by NSF Grant \#DMS-9988809}
\author[A. Wilkinson]{Amie
Wilkinson}
\address{Mathematics Department
  \\ Northwestern University \\
Evanston Illinois, 60208, U.S.A.}
 \thanks{Wilkinson's work was partly funded by NSF Grant \#DMS-0100314}
\email{wilkinso@math.northwestern.edu}
\title{Partial Differentiability of
Invariant Splittings}

\begin{abstract}
A key feature of a general nonlinear
partially hyperbolic  dynamical
system is the   absence of
differentiability of its invariant
splitting.  In this paper, we show that
often partial derivatives of the
splitting exist and
   the splitting depends smoothly
on the dynamical system itself.
\end{abstract}

\thispagestyle{empty}
\maketitle

\bigskip
\bigskip
\bigskip
\bigskip

\centerline{\em Dedicated to David Ruelle on his 65th birthday.}

\thispagestyle{empty}
\vfill
\begin{center}
\today
\end{center}

%%%%%%%%%

\newpage

\section{Introduction}

One of the major technical barriers to the understanding of Anosov
diffeomorphisms is the fact that unstable bundles are not in
general differentiable along stable bundles. This situation persists
for partially hyperbolic diffeomorphisms, where there are also
center bundles present. Under mild bunching conditions, however,
the unstable bundles are differentiable along the center bundles,
see Theorem A below. This fact has already been
observed and exploited in several special situations. First, for
Anosov diffeomorphisms themselves, the unstable bundles are
differentiable with respect the diffeomorphism, as long as partial
derivatives are taken in certain dynamically defined directions
given by conjugating maps \cite{M}. Consequently entropy and SRB states
also
vary differentiably with parameters for Anosov diffeomorphisms and flows
\cite{KKPW,LMM,ruelle1}.  Differentiability of the unstable bundle along
the
center was a crucial ingredient in proving stable ergodicity for
many partially hyperbolic diffeomorphisms \cite{GPS, W, PS, PS2}.  It was
also an ingredient in the construction of nonuniformly hyperbolic
diffeomorphisms
with pathological foliations \cite{SW, D, ruelle2}.  While we have
similar applications in mind for Theorem A, we
will content ourselves here with some general theorems.  We describe the
main results of this paper in the following section; the proofs
occupy the remaining sections.

\section{Statements of Results}\label{s.results}

Suppose  that $f : M \rightarrow  M$ is a partially hyperbolic
diffeomorphism.  The tangent bundle splits as $E^u \oplus E^c
\oplus E^s$.  In general the splitting is continuous but not
$C^1$.   Here we show that under some mild pointwise bunching
conditions, $E^u$ is continuously differentiable in the $E^{c}$
direction, i.e.,
$$ \frac{\partial E^u(p)}{\partial E^c}
$$
exists and is a continuous function of $p \in  M$.   More precisely,
we prove:
\begin{ThmA}
Suppose $f:M\to M$ is $C^2$ and partially hyperbolic with
splitting $TM=E^u\oplus E^c\oplus E^s$. Then, under the
pointwise bunching condition
\begin{equation}
\begin{split}
\label{e:spectralu}
  \sup_p
\frac{\norm{T^c_pf}}
{\boldsymbol{m}(T^u_pf)
\boldsymbol{m}(T^c_pf)}  <  1,
\end{split}
\end{equation}
$E^u$ is continuously differentiable with respect to $E^c$.
\end{ThmA}

Theorem A is a corollary of a more general result about partial
differentiability of dominated splittings -- see Theorem~\ref{t.bundles}
in Section~\ref{s.bundles}.

Next we show, under the same  bunching hypothesis,  that in a family
$t\mapsto f_t$ of
partially hyperbolic diffeomorphisms, the unstable bundle
$E^u(f_t)$
is always continuously differentiable along ``dynamically defined''
curves in $M$.
Roughly speaking, a dynamically defined curve $t\mapsto\varphi_p(t)$
through
$p\in M$ is a $C^1$ curve along which the hyperbolic component of the
dynamics of $f_t$ varies as little as possible.  For example, if $f_t$ is
Anosov and $h_t:M\to M$ is the conjugacy from $f_t$ to $f_0$,
so that $h_t f_0 = f_t h_t$, then $t\mapsto h_t(p)$ is a dynamically
defined
curve.  In the language of Section~\ref{s.functions},  a dynamically
defined
curve is the $M$-component of an integral curve of the center distribution
${\mathbb E}^c$ of the
evaluation map $Eval: M\times I \to M\times I$:
$$(p,t) \mapsto (f_tp, t).$$
We prove that dynamically defined curves always exist, and unstable
bundles,
subject to a bunching condition, vary in a $C^1$ way along them.

\begin{ThmB} Let $\{f_t:M\to M\}_{t\in(-\eps,\eps)}$ be a $C^2$ family of
$C^2$,  partially hyperbolic diffeomorphisms having,
for each $t\in(-\eps,\eps)$, a $Tf_t$-invariant splitting:
$$TM = E^u(f_t)\oplus E^c(f_t) \oplus E^s(f_t).$$

Then there exists $\epsilon_0>0$ so that, for every $p\in M$
there exists a $C^1$ path
$$\varphi_{p}:(-\epsilon_0,\epsilon_0)\to M$$
with $\varphi_p(0) = p$, and with the following property.

If the pointwise bunching condition
\begin{equation}
\begin{split}
  \sup_p
\frac{\norm{T^c_pf_0}} {\boldsymbol{m}(T^u_pf_0)
\boldsymbol{m}(T^c_pf_0)}  <  1
\end{split}
\end{equation}
holds, then $t\mapsto E^u_{\varphi_{p}(t)} (f_t)$ is $C^1$.
\end{ThmB}

Theorem B  follows from a more general result, Theorem~\ref{t.functions},
which states that any invariant, dominated subbundle of a partially
hyperbolic diffeomorphism
is  continuously differentiable along dynamically defined paths,
subject to a bunching condition on the bundle.
In addition, Theorem~\ref{t.functions}  produces,
for any $v\in E^c_p(f_0)$, a dynamically defined path $\varphi_{p,v}$
so  that  $\dot\varphi_{p,v}(0) \in v + E^u\oplus E^s(f_0)$.

The machinery behind the proofs of Theorems A and B is
Theorem~\ref{t.sections}, a refinement of  the $C^1$ Section Theorem
from \cite{HPS}
that  handles  partial derivatives of a section.

In Section~\ref{s.Euisdiff}, we address the question of when
$t\mapsto E^u_p(f_t)$ is differentiable at $t=0$.  The issue here
is of a slightly different nature than that in Theorems A and B.
While $t\mapsto E^u(f_t)$ is always continuously differentiable along
dynamically defined paths, the requirement that the
constant path $t\mapsto p$ be dynamically defined for all $p$ is
a stringent one,  satisfied only for very special families.

If, instead of requiring that $t\mapsto E^u_p(f_t)$  be $C^1$ in a
given family,
we just ask that it be differentiable at $t=0$ {\em but for all families
through $f_0$},
then  the actual dynamics of $f_0$ becomes irrelevant.
It is easy to see that $p\mapsto E^u_p(f_0)$ must be $C^1$ for
this property to hold.  What is interesting is that nonsmoothness
of  $p\mapsto E^u_p(f_0)$ is the {\em only} obstruction to  the
differentiability of $t\mapsto  E^u_p(f_t)$ at $t=0$ in every family.
Building on Theorem B, one can show:

\begin{ThmC} Let $\{f_t:M\to M\}_{t\in(-\eps,\eps)}$ be a $C^1$ family of
$C^2$,  partially hyperbolic diffeomorphisms having,
for each $t\in(-\eps,\eps)$, a $Tf_t$-invariant splitting:
$$TM = E^u(f_t)\oplus E^c(f_t) \oplus E^s(f_t).$$

Assume that the pointwise bunching condition
\begin{equation}
\begin{split}
  \sup_p
\frac{\norm{T^c_pf_0}} {\boldsymbol{m}(T^u_pf_0)
\boldsymbol{m}(T^c_pf_0)}  <  1
\end{split}
\end{equation}
holds.  Assume also that $E^u(f_0)$ is a $C^{2-\eps}$ subbundle of $TM$,
for all $\eps>0$.

Then  for all $p\in M$, $$t\mapsto E^u_p(f_t)$$
  is differentiable at $t=0$.

If  $\varphi_p$ is any dynamically defined path through $p\in M$
given by Theorem B,  then:
$$E^u_p(f_t) - E^u_p(f_0) =
\left(\frac{d}{dt}E^u_{\varphi_{p}t}(f_t)\,\vert_{t=0}-D_pE^u(f_0)(\frac{d}{dt}
\varphi_{p}t\,\vert_{t=0})\right)t + O(t^{1 + \eta}),$$ for some
$\eta>0$.
\end{ThmC}

Subsequent to proving Theorem C, we learned of a more general result,
due to Dolgopyat:

\begin{ThmD}[Dolgopyat, \cite{D}, Theorem 3] Let $\{f_t:M\to
M\}_{t\in(-\eps,\eps)}$ be a $C^1$ family of
$C^2$ diffeomorphisms having,
for each $t\in(-\eps,\eps)$, a $Tf_t$-invariant dominated splitting:
$$TM = R(f_t)\oplus S(f_t) \oplus T(f_t).$$
Suppose that $p\mapsto S_p(f_0)$ is $C^1$.  Then,
for every $p\in M$, $t\mapsto S_p(f_t)$ is differentiable at $t=0$.
\end{ThmD}

Dominated splittings are defined in Section~\ref{s.bundles}.  In
particular,
Theorem D applies when
$S$ is $E^u$, $E^c$, $E^s$,  $E^{cu}$, or $E^{cs}$.
We present an exposition of Dolgopyat's proof of Theorem D in
Section~\ref{s.Euisdiff}.

\section{Partial Derivatives of an
Invariant Section} \label{s.sections} Let
$$
\begin{CD}
$$
  V @>\text{\normalsize $\qquad
F \qquad$}>> V\\
  @V\text{\normalsize $\pi$}VV
@VV\text{\normalsize
$\pi$}V\\
  M @>\text{\normalsize $\qquad
f \qquad$}>> M$$
\end{CD}
$$
be a $C^1$ fiber preserving  map,
where
$V$ is a smooth, finite dimensional
fiber bundle over the compact
manifold
$M$, and $f$ is a diffeomorphism.
In addition   assume that there is
a section $\sigma  : M \rightarrow
V$, invariant under $F$ in the sense
that
$$
F(\sigma (p)) = \sigma (f(p))
$$
for all $p \in  M$.

In general there is no reason that
$\sigma $ is smooth, or even
continuous.  For example, if $F$ is
the identity map, every section of
$V$ is $F$-invariant.  In
\cite{HPS}, we showed that if $V$ is
a Banach bundle and
$F$ is a fiber contraction then
$\sigma $ is unique and continuous,
and furthermore, if the fiber
contraction dominates the base
contraction sufficiently then the
$\sigma $ is of class $C^r$.

Since $F$ preserves fibers, $TF$
preserves the
``vertical'' subbundle,
$
\operatorname{Vert} \subset  TV
$
whose fiber at $v \in  V$ is
$\operatorname{kernel}T_v\pi$.  We
write
$
T^{\operatorname{Vert}}_vF$ for
the restriction of $T_vF$ to
$\operatorname{Vert}_v$,
$$
T^{\operatorname{Vert}}_v F :
\operatorname{Vert}_v \rightarrow
\operatorname{Vert}_{Fv}.
$$

We assume that $TV$ carries a
Finsler structure and that
$
k_p =
\norm{T^{\operatorname{Vert}}
_{\sigma p}F}$
has
$$\sup_{p \in  M} k_p < 1,
$$
which means that $F$ is a fiber
contraction in the neighborhood of
$\sigma M$.

%%--- THEOREM ----
\begin{Thm}
\label{t.sections} Suppose   that $E \subset  TM$ is a continuous
$Tf$-invariant subbundle such that
$$
\sup_{p \in  M} k_p
\norm{(T_p^Ef)^{-1}} < 1
$$
where
$T^Ef$ is the restriction of $Tf$ to
$E$.
Then
$\sigma $  is
continuously differentiable in the
$E$-direction in the sense that there
is a continuous   map $ H : E
\rightarrow  TV$ such that
\begin{itemize}
   \item[(a)]

For each $p  \in  M$,  $H : E_p
\rightarrow T_{\sigma p}V$ is linear.

\item[(b)]

$T\pi \circ  H =
\operatorname{Id} : E
\rightarrow  E$.
   \item[(c)]

If $\gamma $ is a $C^1$
arc in $M$ that is everywhere
tangent to $E$ then
$$
(\sigma \circ \gamma )^{\prime}(t)
= H(\gamma ^{\prime}(t)).
$$
\end{itemize}
In particular, if $E$ is integrable
then the restriction of $\sigma $ to
each
$E$-leaf is
$C^1$.
\end{Thm}

We refer to $H$ as the partial
derivative of $\sigma $ in the
$E$-direction
$$
H = \frac{\partial \sigma }{\partial
E}.
$$

\begin{Rmk} If, in addition, there exist $C^r$ submanifolds everywhere
tangent
to $E$, for some $r\in(0,\infty)$, then  $C^r$ smoothness of $\sigma
$ along $E$ (i.e., along
these manifolds) can be assured by assuming that
$$
\sup_{p}
k_p\norm{(T^E_pf)^{-1}}^r < 1.
$$
\end{Rmk}

\begin{Rmk}
When $E$ is integrable, the proof of Theorem~\ref{t.sections} is a
fairly simple application of the Invariant Section Theorem of
\cite{HPS}.
  It is the
non-integrable case that requires
some new ideas.
\end{Rmk}

\begin{Rmk}
There is a uniformity about
$\partial \sigma /\partial E$.  (In the
integrable case, this uniformity is
automatic.)   Fix
$p
\in  M$ and extend each
$w
\in  E_p$ with $\abs{w} \leq  1$ to  a
continuous vector field $X_w$
everywhere subordinate to $E$,  and
do so in a way that depends
continuously on $w$. Let
$\gamma _w$ be an integral curve of
$X_w$ through $p$.  Since $E$ is
only continuous, the integral curve
$\gamma _w$ need not be uniquely
determined by
$X_w$.   Nevertheless, for all $p$ in
any fixed
$C^1$ chart, as $t \rightarrow  0$ we
have
$$
\frac{\sigma \circ \gamma _w(t) -
\sigma p}{t}\rightarrow
H(w)
$$ uniformly.
\end{Rmk}

\begin{Rmk}
Since $M$  is finite dimensional,
Peano's Existence Theorem implies
that there exist $C^1$ arcs
  everywhere tangent to a continuous
plane field, and thus the hypothesis of assertion (c) in
Theorem~\ref{t.sections} is satisfied.  In the infinite dimensional
case, however, Peano's Theorem fails and (c) could become vacuous.
\end{Rmk}

%%---  Proof --
\begin{proof}
[Proof of Theorem~\ref{t.sections}] We proceed by the graph
transform techniques in \cite{HPS}. Choose a continuous subbundle
$\operatorname{Hor} \subset TV$, complementary to
$\operatorname{Vert}$,
$$
\operatorname{Hor} \oplus
\operatorname{Vert} = TV.
$$
For example, we could introduce a
Riemann structure on $TV$ and take
$\operatorname{Hor}_v$ as the
orthogonal complement to
$\operatorname{Vert}_v$.
Note that $T\pi $ sends each
subspace
$\operatorname{Hor}_v$
isomorphically onto $T_pM$, $p =
\pi v$.  With respect to the horizontal /
vertical splitting we write
\begin{equation*}
\begin{split}
T_{v}F =
\begin{bmatrix}
A_v  & 0  \\
C_v  & K_v
\end{bmatrix} =
\begin{bmatrix}
A_v : \operatorname{Hor} _v
\rightarrow
\operatorname{Hor} _{Fv}  & 0
\\
C_{v} : \operatorname{Hor}
_v \rightarrow
\operatorname{Vert} _{Fv}  & K_v :
\operatorname{Vert} _v
\rightarrow \operatorname{Vert}
_{Fv}
\end{bmatrix}.
\end{split}
\end{equation*}

Let $\overline{E} \subset
\operatorname{Hor}$ be the lift of
$E$.  That is, $T\pi $ sends the plane
$\overline{E}_v$ isomorphically to
$E_p$,  $p = \pi v$.  Since $E$ is
$Tf$-invariant and $F$ covers $f$,
$\overline{E}$ is $A$-invariant in
the sense that
$$
\begin{CD}
$$
  \overline{E}_v @>\text{\normalsize
$\qquad A_v \qquad$}>>
\overline{E}_{Fv}\\
  @V\text{\normalsize $T\pi $}VV
@VV\text{\normalsize
$T\pi $}V\\
  E_p @>\text{\normalsize $\qquad
Tf \qquad$}>> E_{fp}
$$
\end{CD}
$$
commutes.

Let $L$ be the bundle over $M$
whose fiber at $p$ is
$$
L_p = L(\overline{E}_{\sigma p},
\operatorname{Vert}_{\sigma p}).
$$
An element in $ L_p$ is a linear
transformation $P :
\overline{E}_{\sigma p} \rightarrow
\operatorname{Vert}_{\sigma p}$.
Let
$LF$ be the \textbf{graph transform}
on
$L$ that sends $P \in
L_p$ to
$$
P^{\prime} = (C_{\sigma p} +
K_{\sigma p}P)
(A_{\sigma
p}|_{\overline{E}_{\sigma p}})^{-1}
\in L_{\sigma p}.
$$
Then  $TF$ sends the graph of
$P$ to the graph of
$P^{\prime}$ and $LF$ is an affine
fiber contraction
$$
\begin{CD}
$$
  L @>\text{\normalsize $\qquad
LF \qquad$}>> L
\\
  @V\text{\normalsize $\pi $}VV
@VV\text{\normalsize
$\pi $}V\\
  M @>\text{\normalsize $\qquad
f \qquad$}>> M.
$$
\end{CD}.
$$
By \cite{HPS}, $L$ has a unique
$LF$-invariant section $\Lambda  :
M
\rightarrow L$, and $\Lambda $ is
continuous.   Define $H_p : E_p
\rightarrow T_{\sigma p}V$ by
commutativity of
$$
\begin{CD}
$$
  \overline{E}_{\sigma p}
@>\text{\normalsize
$\qquad \operatorname{Id}_p \oplus
\Lambda _p
\qquad$}>> \overline{E}_{\sigma
p} \oplus
\operatorname{Vert}_{\sigma p}\\
  @V\text{\normalsize $T\pi $}VV
@VV\text{\normalsize
$\operatorname{Inclusion}$}V\\
  E_p @>\text{\normalsize $\qquad
H_p \qquad$}>> T_{\sigma p}V
$$
\end{CD}
$$
where $\operatorname{Id}_p$ is the
identity map $\overline{E}_{\sigma
p} \rightarrow  \overline{E}_{\sigma
p}$.  Then $H : E \rightarrow  TV$
is   the unique bundle map such that
$HE$ is a
$TF$-invariant subbundle of
$T_{\sigma M}V$,
$$
\begin{CD}
$$
  HE @>\text{\normalsize $\qquad
TF \qquad$}>> HE\\
  @V\text{\normalsize $T\pi $}VV
@VV\text{\normalsize
$T\pi $}V\\
  E @>\text{\normalsize $\qquad
Tf \qquad$}>> E
$$
\end{CD}
$$
commutes,  and $T\pi  \circ  H =
\operatorname{Id}_E$.
  We claim that $H$  is
the partial derivative of
$\sigma $ in the $E$-direction.

Let $\gamma  : (a, b)
\rightarrow  M$ be a
$C^1$ arc such that $\gamma
$ is everywhere tangent to $E$.
   To complete
the proof of the theorem, we must
show that
$$
(\sigma \circ \gamma )^{\prime}(t)
= H(\gamma ^{\prime}(t)).
$$
For $n \in  \mathbb{Z}$, set
$\gamma _n =f^n \circ
\gamma
$ and
$$
\Gamma = \bigsqcup_{n
\in  \mathbb{Z}}
\gamma _n.
$$
   This means that we
consider the disjoint union of the
arcs $\gamma _n$, so  if two of them
cross in $M$, we ignore the
crossing in $\Gamma $.  The one
dimensional manifold
$\Gamma
$ is   noncompact; it has countably
many components  $\gamma _n$.  In
the same way, we discretize
  $V$ as
\begin{equation*}
\begin{split}
V_{\Gamma } = \bigsqcup_n
V|_{\gamma _n}.
\end{split}
\end{equation*}
We equip $V_{\Gamma }$ and
$T\Gamma $ with the Finslers they
inherit from $V$ and $M$.  Then
$F_{\Gamma } = F|_{V_{\Gamma
}}$ is a fiber contraction
$$
\begin{CD}
$$
  V_{\Gamma } @>\text{\normalsize
$\qquad F_{\Gamma } \qquad$}>>
V_{\Gamma }\\
  @V\text{\normalsize $\pi $}VV
@VV\text{\normalsize
$\pi $}V\\
  \Gamma  @>\text{\normalsize
$\qquad f
\qquad$}>> \Gamma ,
$$
\end{CD}
$$
and the fiber contraction dominates
the base contraction since
$$
\sup_pk_p\norm{T^E_pf} < 1
$$
and $T\Gamma  \subset  E$.
Furthermore, $F_{\Gamma }$ is
uniformly $C^1$ bounded since
$M$ is compact.  The Invariant
Section Theorem of \cite{HPS} then
implies that $V_{\Gamma }$ has a
unique $F_{\Gamma }$-invariant
section $\sigma _{\Gamma }$, and
$\sigma _{\Gamma }$ is of class
$C^1$.  Furthermore the tangent
bundle of $\sigma _{\Gamma
}(\Gamma )$ is the unique nowhere
vertical
$TF_{\Gamma }$-invariant line
field in $TV_{\Gamma }$.

The restriction of $\sigma $
to $\Gamma = \bigsqcup_n\gamma
_n$ is $F_{\Gamma }$-invariant, so
  by uniqueness
$$
\sigma _{\Gamma } = \bigsqcup_n
\sigma |_{\gamma _n}.
$$
We claim that
$$
H(T\Gamma ) =  T(\sigma
_{\Gamma }\Gamma ).
$$
Again the reason is uniqueness. We
know that
$T(\sigma _{\Gamma }\Gamma )$ is
the unique $TF_{\Gamma
}$-invariant, nowhere vertical line
field  defined over $\sigma
_{\Gamma }\Gamma
$.  But   commutativity of
$$
\begin{CD}
$$
  HE
@>\text{\normalsize
$\qquad TF_{\Gamma } \qquad$}>>
HE\\
  @A\text{\normalsize $H $}AA
@AA\text{\normalsize
$H $}A\\
  E @>\text{\normalsize $\qquad
Tf \qquad$}>> E
\\
  @A\text{\normalsize
$\operatorname{Inclusion}$}AA
@AA\text{\normalsize
$\operatorname{Inclusion}$}A\\
  T\Gamma  @>\text{\normalsize
$\qquad Tf \qquad$}>> T\Gamma
$$
\end{CD}
$$
implies that
$H(T\Gamma )$ is a second such
line field.  By uniqueness they are
equal.

To complete the proof, we show that
the line field equality implies
the vector equality
$$
\frac{d}{dt} \; \sigma \circ \gamma
(t) = H(\gamma ^{\prime}(t)),
$$
as the theorem asserts.
Differentiating
$\gamma (t) = \pi  \circ  \sigma
_{\Gamma } \circ \gamma (t)$ gives
$$
\gamma ^{\prime}(t) = T \pi  \circ
T  \sigma _{\Gamma }
(\gamma ^{\prime}(t)).
$$
The vector $T  \sigma _{\Gamma }
(\gamma ^{\prime}(t))$ lies in the
span of
$H(\gamma ^{\prime}(t))$, so there
is a real number $c(t)$ such that
$T  \sigma _{\Gamma } (\gamma
^{\prime}(t))  = H( c(t) \gamma
^{\prime}(t))$.
This gives
$$
\gamma ^{\prime}(t) = T\pi \circ
H(c(t) \gamma ^{\prime}(t)).
$$
Since $T\pi  \circ  H =
\operatorname{Id}_E$ we have
$$
\gamma ^{\prime}(t) = c(t) \gamma
^{\prime}(t)
$$
and $c(t) = 1$.  Thus
$$
\frac{d}{dt} \; \sigma _{\Gamma }
\circ \gamma (t)
=
T  \sigma _{\Gamma } (\gamma
^{\prime}(t))
=
H(c(t)\gamma ^{\prime}(t)) =
H(\gamma ^{\prime}(t)).
$$
\end{proof}
%%  --- QED ---

\begin{Rmk}
Above, it is
assumed that $\gamma $ is everywhere
tangent to $E$.  One might expect that
tangency of $\gamma $ to $E$ at $p
= \gamma (0)$ suffices to prove
that $(\sigma  \circ  \gamma
)^{\prime}(0) = H(\gamma
^{\prime}(0))$.  This is not so.  For
example
$E$ can be the flow direction of an
Anosov flow.  The bundle $E^u$ can
be H\"older, but not $C^1$.  Say its
H\"older exponent is
$\theta < 1$.  One can construct a
$C^1$ curve  $\gamma (t)$
which is tangent to
$E$ at $p = \gamma (0)$, but which
diverges from
$E$ at a rate $t^{1+ \epsilon }$.  The
difference between
$E^u_{\gamma (t)}$  and
$E^u_p$ is then on the order of
$  t^{\theta +\epsilon \theta }$.  If
$\epsilon $ is small this exponent  is
$< 1$, and the map $t \mapsto
E^u_{\gamma (t)}$ fails to be
differentiable at $t = 0$.
\end{Rmk}

\section{A Series Expression for
$\partial \sigma / \partial E$}
\label{s:formula}
As
above $\sigma $ is the unique
$F$-invariant section and  $H =
\operatorname{Id}_E \oplus
\Lambda $ is its partial derivative in
the
$E$-direction.   Naturally,  $\partial
\sigma  / \partial E$  depends on the
choice of horizontal subbundle
$\operatorname{Hor} \subset TV
$.  We use the isomorphism
$T\pi : \operatorname{Hor}_{\sigma
p}
\rightarrow
  T_pM$ to identify the linear
map
$A_{\sigma p} :
\operatorname{Hor}_{\sigma
p}\rightarrow
\operatorname{Hor} _{\sigma (fp)}$
with its
$T\pi$-conjugate  $T_pf$.  Then,
using the
canonical isomorphism
$\operatorname{Vert} _{\sigma p}
\approx V_p$, we can express $TF =
\begin{bmatrix} A & 0 \\ C &K
\end{bmatrix}$ as
\begin{equation*}
\begin{split}
T_{\sigma p}F =
\begin{bmatrix}
T_pf : T_pM \rightarrow  T_{fp}M
& 0
\\
C_p : T_pM \rightarrow V_{fp}
   & K_p : V_p
\rightarrow V_{fp}
\end{bmatrix}.
\end{split}
\end{equation*}
Thus, the bundle map $LF : L
\rightarrow  L$ becomes
$$
P \mapsto (C_p + K_pP) \circ
(T^E_{fp}f^{-1}).
$$
Denote by $\Lambda _0$ the zero
section of $L$, and call its
$N^{\textrm{th}}$ iterate in $L$,
$$
\Lambda  _N =
(LF)^N(\Lambda  _0).
$$
We know that $\Lambda  _N
\rightarrow \Lambda  $ uniformly as $N
\rightarrow \infty$.  Also,  we claim
that
$$
\Lambda  _N(p) =
\sum_{n=0}^{N-1} K_p^{n} \circ
C_{f^{-n-1}(p)} \circ
(T_p^Ef^{-n-1})
$$
where $K^0 = \operatorname{Id}$
and for $n  \geq 1$,
$$
K_p^{n} = K_{f^{-1}(p)}
\circ  \dots \circ K_{f^{-n}(p)}
: V_{f^{-n}(p)} \rightarrow
V_p.
$$
If $N =
1$ we have
$$
\Lambda  _1(p) = (C_{f^{-1}(p)} +
K_{f^{-1}(p)}P_0)(T_p^Ef^{-1}) =
C_{f^{-1}(p)} T_p^Ef^{-1}
$$
because $\Lambda _0 = 0$
implies that $P_0 = 0$.  Thus, the
assertion holds with
$N = 1$;  the proof is completed by
induction.

Since the partial sums of the infinite
series
$\sum_{n=0}^{\infty} K_{p}^{n}
C_{f^{-n-1}(p)} T_p^Ef^{-n-1}$
converge uniformly to $\Lambda  $,
we are justified in writing
$$
  \frac{\partial \sigma
}{\partial  E} =  H(p) =
\operatorname{Id}_E \oplus
\sum_{n=0}^{\infty} K_{p}^{n}
C_{f^{-n-1}(p)} T_p^Ef^{-n-1}.
$$

\section{The Dominated and Partially Hyperbolic
Cases: Proof of Theorem A}\label{s.bundles}

Let $f:M\to M$ be a diffeomorphism of a compact Riemannian (or
Finslered) manifold. In this section, we consider $Tf$-invariant
dominated and partially hyperbolic splittings of the tangent
bundle $TM$. We recall the definitions. The \textbf{conorm}  of a
linear transformation $T : X \rightarrow  Y$ is
$$
\boldsymbol{m}(T) = \inf_x \frac{\abs{Tx}}{\abs{x}}.
$$
Suppose that $TM$ splits
as a $Tf$-invariant sum of two bundles:
$$TM = R\oplus S.$$
This splitting is {\em dominated} if the following condition
holds:
\begin{equation*}
\begin{split}
\inf_p \frac{\boldsymbol{m} (T_p^Rf)}{\norm{ T^S_pf} }> 1
\end{split}
\end{equation*}
where the notation $T^Xf$ is used for the restriction of $Tf$ to
the bundle $X$.  We say that $f$ has a {\em dominated
decomposition} if there is a $Tf$-invariant dominated splitting of
$TM$.

Then $f$ is {\em partially hyperbolic} if it has a $Tf$-invariant splitting
$TM=E^u\oplus E^c\oplus E^s$ such that:
\begin{enumerate}
\item $E^u\oplus ( E^c\oplus E^s) $ and $(E^u\oplus E^c)\oplus E^s$
are both dominated splittings of $TM$,
\item \begin{equation*}
\begin{split}
\inf_p\boldsymbol{m} (T^u_pf) > 1\qquad & \sup_p\norm{(T^s_pf)} <
1
\end{split}
\end{equation*}
\end{enumerate}
In other words,  $Tf$ expands   vectors in $E^u$, contracts
vectors in $E^s$, and is relatively neutral on vectors in $E^c$.

Theorem~\ref{t.sections} can be applied in the dominated
decomposition and partially hyperbolic contexts, as follows.

\begin{Thm}\label{t.bundles}
Suppose that the $C^2$ diffeomorphism $f:M\to M$ has a dominated
decomposition:
$$TM = R\oplus S,$$
And let $E\subset TM$ be a  $Tf$-invariant subbundle.
Then, under the pointwise bunching condition
\begin{equation}
\begin{split}
\label{e:spectral}
  \sup_p
\frac{\norm{T^S_pf}} {\boldsymbol{m}(T^R_pf)
\boldsymbol{m}(T^E_pf)}  <  1,
\end{split}
\end{equation}
$R$ is continuously differentiable with respect to $E$.
\end{Thm}

Theorem A is an immediate corollary of Theorem~\ref{t.bundles} ,
where we set $R=E^u$, $S=E^c\oplus E^s$, and $E=E^c$.

  In \cite{HPS} we defined the ``bolicity'' of a linear transformation
to be the ratio of its norm to its conorm, so the poinwise bunching
condition in Theorem A can be re-stated as
$$
\sup_{p} \frac{\operatorname{bol}(T^c_pf)}
{\boldsymbol{m}(T_p^uf)} < 1.
$$
  Similarly, to show
that $E^s$ is differentiable along $E^c$, we assume
\begin{equation}
\begin{split}
\label{e:spectrals}
  \sup_p
\norm{T^s_pf} \operatorname{bol}(T^c_pf)   < 1.
\end{split}
\end{equation}

\bigskip

\begin{proof}
[Proof of Theorem~\ref{t.bundles}] Let $d$ be the fiber dimension
of $R$ in the dominated decomposition
$$TM=R\oplus S.$$
$Tf$  acts naturally on the Grassmann $G = G(d, TM)$ of all
$d$-planes in $TM$,
$$
Gf  :  G  \rightarrow  G.
$$
If $\Pi$ is a $d$-plane in $T_pM$ then $Gf(\Pi) = T_pf(\Pi)$.
  Since
$f$ is
$C^2$,
$Gf$  is a
$C^1$ fiber preserving  map,
$$
\begin{CD}
$$
  G @>\text{\normalsize $\qquad
Gf \qquad$}>> G\\
  @V\text{\normalsize $\pi$}VV
@VV\text{\normalsize
$\pi$}V\\
  M @>\text{\normalsize $\qquad
f \qquad$}>> M
$$
\end{CD},
$$
where $\pi$ sends  $\Pi \subset T_pM$ to $p$. Besides, $p \mapsto
R_p$ is a $Gf$-invariant section of $G$.  We show that, at the
invariant section $R$, the fiber contraction rate dominates the
contraction rate along $E$ as follows.

A compact neighborhood $N_p$ of $R_p$ in $G_p$ consists of
$d$-planes $\Pi $ such that $\Pi = \operatorname{graph}P$ where $P
: R_p \rightarrow  S_p$ is a linear transformation with
$\norm{P} \leq  1$.  Give $G$ a Finsler which is the operator norm
on each $N_p$ and any other Finsler on the rest of $G$.  Then $Gf$
is a fiber preserving map whose fiber contraction rate at $R_p$ is
$$
k_p = \frac{\norm{T^{S}_pf}} {\boldsymbol{m}(T^R_pf)} < 1.
$$
Since this dominates the contraction rate $\boldsymbol{m}(T^E_pf)$
along $E$,  Theorem~\ref{t.sections} applies and $p \mapsto R_p$
is seen to be a continuously differentiable function of $p$ in the
$E$ direction.
\end{proof}

%%%%%%%%%%%%%%%%%%%%%%%%%%%%%%%%%%%%%%%%%%%%%%%%%%%%%%%%%%%%%%%%%%%5

\section{A series formula for $\partial E^u/ \partial
E^c$}\label{s.phformula}
  From Section~\ref{s:formula} we know that
there is a series that expresses $\partial E^u/ \partial E^c$.  We
write this formula out after making   a convenient choice of the
horizontal bundle.

To do so, we  coordinatize $G$  near
$E^u_p$ as follows.
Fix a smooth Riemann structure on
$M$ that exhibits the partial
hyperbolicity of $f$, and let
$\operatorname{exp}$ be its
exponential map.  Abusing notation,
we denote by
$\mathbb{R}^u$ and
$\mathbb{R}^{cs}$ the planes
$\mathbb{R}^u \times  0$ and $0
\times  \mathbb{R}^{cs}$ in
$\mathbb{R}^m$.  For each
$p
\in  M$, define a linear map $I_p :
\mathbb{R}^m \rightarrow  T_pM$
that carries $\mathbb{R}^u $
and $  \mathbb{R}^{cs}$
isometrically to $E^u_p$ and
$E_p^{cs}$.  The restriction of
$\operatorname{exp}_p \circ  I_p$
to a small neighborhood $U  =U_p$
of
$0$ in $\mathbb{R}^m$ is a
diffeomorphism of $U$ to a
neighborhood $Q = Q_p$ of $p$ in
$M$,
$$
\varphi_p : U \rightarrow  Q
$$
  and
\begin{itemize}
   \item[(a)] $\varphi_p(0) = p$
   \item[(b)] $T_0\varphi_p$ carries
$ \mathbb{R}^u
  $ and $
\mathbb{R}^{cs}$ isometrically to
$E^u_p$ and
$E_p^{cs}$.
\item[(c)]  If we denote
by $E_{pq}^u$ and
$E_{pq}^{cs}$ the
planes
$T_x\varphi_p(\mathbb{R}^u)$ and
$T_x\varphi_p(\mathbb{R}^{cs})$,
where $q = \varphi_p(x)$,  then $q
\mapsto E_{pq}^{u}
\oplus E_{pq}^{cs}$ is a smooth
splitting of $TQ$ that
reduces to $E^u_p \oplus E_p^{cs}$
when $p = q$.
\end{itemize}

We now  coordinatize $G$ near
$E^u_p$.  Let
$\mathcal{M}$ be the space of
$(u \times cs)$-matrices, thought of
as linear transformations
$X : \mathbb{R}^u \rightarrow
\mathbb{R}^{cs}$.  Given
$(x, X)
\in U
\times
\mathcal{M}$, let $q = \varphi
_p(x)$ and consider the linear
transformation $S : E^u_{pq}
\rightarrow E^{cs}_{pq}$,  defined
by commutativity of
$$
\begin{CD}
$$
\mathbb{R}^u  @>\text{\normalsize
$\qquad X \qquad$}>>
\mathbb{R}^{cs}
\\
  @V\text{\normalsize $T_x\varphi
_p$}VV @VV\text{\normalsize
$T_x\varphi  _p$}V\\
  E^u_{pq}
@>\text{\normalsize
$\qquad  S
\qquad$}>>E^{cs}_{pq} .
$$
\end{CD}
$$
The graph of $S$ is a plane   $\Pi
\in G$ near $E^u_p$, and thus
$$
\Phi _p : (x, X) \mapsto  \Pi
$$
is a local trivialization of $G$ at
$E^u_p$.

Because   $U \times
\mathcal{M}$ is a product, $T(U
\times  \mathcal{M})$  carries a
natural horizontal structure,   the
horizontal space at
$(x, X)$  being
$$
\mathbb{R}^m \times  0 \; \subset \;
\mathbb{R}^m \times \mathcal{M}
= T_{(x, X)}
(U \times \mathcal{M}).
$$
We define the horizontal space at
$\Pi = \Phi _p(0, X)  \in G_p$ to be
$$
\operatorname{Hor}_{\Pi }
=
T_{(0, X)}\Phi _p(\mathbb{R}^m
\times  0).
$$

Writing $T(Gf) : TG \rightarrow TG$ with respect to the horizontal
/ vertical splitting of $TG$
  gives
$$
T_{\Pi }(Gf) =
\begin{bmatrix}
  A_{\Pi } :
\operatorname{Hor}_{\Pi }
\rightarrow
\operatorname{Hor}_{Gf(\Pi )}  &
0  \\
  C_{\Pi } : \operatorname{Hor}_{\Pi
} \rightarrow
\operatorname{Vert}_{Gf(\Pi )}
&  K_{\Pi } :
\operatorname{Vert}_{\Pi }
\rightarrow
\operatorname{Vert}_{Gf(\Pi )}
\end{bmatrix}
$$
Take $\Pi  = E^u_p$ and
identify
$$
T_{E_p^u}G_p
= \operatorname{Vert}_{E^u_p}
\approx L(E^u_p, E^{cs}_p).
$$
Fix $v \in E^c_p$.  Then
$C_{f^{-n-1}p} \circ
T^cf^{-n-1}(v)$ is a linear
transformation
$Y_n (v): E^u_{f^{-n}p} \rightarrow
E^{cs}_{f^{-n}p}$, and $Y_n(v)$  is
susceptibe to the $n^{\textrm{th}}$
power of the linear graph transform,
which converts it to a linear
transformation
$ E^u_p \rightarrow
E^{cs}_p$ defined by
$$
T^{cs}_{f^{-n}p}f^n
\circ  (Y_n(v)) \circ T_p^uf^{-n}.
$$
This is the same as the repeated action
of $K$.  (That is, the graph
transform of $Tf^n$ is the
same as the $n^{\textrm{th}}$
power of the graph transform
of $Tf$.) Thus, by the formula in
Section~\ref{s:formula},
$$
\frac{\partial E^u}{\partial E^c} (v)
= \sum_{n=0}^{\infty}
T^{cs}_{f^{-n}p}f^n
\circ  (C_{f^{-n-1}p} \circ
T^cf^{-n-1}(v)) \circ T_p^uf^{-n}.
$$

We also can express this
in  charts as follows.
Writing $f$ in the $\varphi $-charts
gives
$$
f_p = \varphi _{fp}^{-1} \circ  f
\circ  \varphi _p
$$
and
$$
(Df_p)_x  =
\begin{bmatrix}
  D^{u,u}_xf_p  &  D^{cs,u}_xf_p
\\
  D^{u,cs}_xf_p  &
D^{cs,cs}_xf_p
\end{bmatrix}
$$
where the $D^{u,u}_xf_p$ block
consists of the partial derivatives of
the
$u$-components of $f_p$ with
respect to the $u$-variables,
evaluated at the point $x$, etc.  At $x
=0$, the off-diagonal blocks are zero,
while the diagonal blocks are
$T\varphi $-conjugate to
$T^u_pf$ and
$T^{cs}_pf$.   Thus, the coordinate
expression of $Gf$ becomes
$$
(Gf)_p : (x, X) \mapsto (f_px,
(D_x^{u, cs}f_p +
(D_x^{cs,cs}f_p)X) (D_x^{u,u}f_p +
(D_x^{cs,u}f_p)X)^{-1}.
$$
Differentiating this
   with respect to
$x$ and
$X$ at the origin $  (0, 0) \in
\mathbb{R}^m \times  \mathcal{M}$
yields
$$
(D((Gf)_p))_{(0,0)} =
\begin{bmatrix}
  A_p : \mathbb{R}^m \rightarrow
\mathbb{R}^m &  0 \\
  C_p : \mathbb{R}^m \rightarrow
\mathcal{M}&   K_p : \mathcal{M}
\rightarrow  \mathcal{M}
\end{bmatrix},
$$
where $A_p$ is $T\varphi
$-conjugate to $T_pf$,
$$
A_p = (T_{0}\varphi _{fp})^{-1}
\circ  T_pf \circ T_0\varphi _p,
$$
$C_p$ represents the second
derivatives of $f$ in the $\varphi
$-charts,
\begin{equation*}
\begin{split}
C_p  &=
\frac{\partial }{\partial
x}\Big|_{x=0} (D_x^{u,
cs}f_p)(D_x^{u,u}f_p)^{-1}
\\
&=
(D_x(D_x^{u,
cs}f_p))(D_x^{u,u}f_p)^{-1}
\\
&\quad
-
(D_x^{u, cs}f_p)
(D_x^{u,u}f_p)^{-1}
(D_x(D_x^{u,u}f_p))
(D_x^{u,u}f_p)^{-1},
\end{split}
\end{equation*}
and, because the off-diagonal blocks
vanish at the origin,
$K_p$ is
$T\Phi
$-conjugate to the graph
transform of
$Tf$,
\begin{equation*}
\begin{split}
P \mapsto T^{cs}_pf \circ P
\circ (T^{u}_pf)^{-1}.
\end{split}
\end{equation*}

It is worth noting that the norm of
$C$ is uniformly bounded on a
neighborhood of $E^u$ in $G$
because $f$ is $C^2$ and $M$ is
compact.  Also, this is clear from the
formula expressing $C$ in the
$\varphi $-charts.

%%%%%%%%%%%%%%%%%%%%%%%%%%

\section{Dependence of $E^u$,
$E^s$ on $f$: Proof of Theorem B} \label{s.functions}

As has been highlighted in the
Katok-Milnor examples \cite{Milnor}, the conjugacy between an
Anosov diffeomorphism and its perturbations is a smooth function
of the perturbation, even though the conjugacies themselves are
only continuous.  For example, consider a $1$-parameter family of
Anosov diffeomorphisms $g_t : M \rightarrow  M$.  The map $g_0$ is
conjugate to $g_t$ by a homeomorphism $h_t : M \rightarrow  M$,
and $h_t$ is uniquely determined by the requirements that $h_0 =
\operatorname{Id}$ and $t \mapsto h_t$ is continuous.  The map
\begin{equation*}
\begin{split}
  (-\epsilon , \epsilon ) \times  M
&\rightarrow  (-\epsilon ,
\epsilon ) \times  M
\\
(t, p) &\mapsto (t, g_tp)
\end{split}
\end{equation*}
is smooth, partially hyperbolic,  and
supports a center foliation
$\mathcal{W}^c$ whose
   leaf through   $(0,p)$ is
$  \{(t, h_t(p)) : t \in (-\epsilon , \epsilon )\}$.  As a
foliation $\mathcal{W}^c$ is only continuous, but its leaves are
smooth. Theorem~\ref{t.sections} applies perfectly well to this
situation, and we conclude that $E^u_{h_tp}$, $E^s_{h_tp}$ are
$C^1$ functions of $t$.  In this section we replace the Anosov
condition by  partial hyperbolicity, and derive an analogous
result.

We assume that $f _0 : M \rightarrow  M$ is a $C^2$, partially
hyperbolic diffeomorphism with splitting
$$TM = E^u\oplus E^c\oplus E^s$$
and that $\mathcal{F}$ is a small neighborhood of $f_0$ in
$\operatorname{Diffeo}^2M$.  By the usual linear graph transform
techniques, all  $f \in \mathcal{F}$ are partially hyperbolic and
their splittings
$$TM = E^u(f)\oplus E^c(f)\oplus E^s(f)$$
depend continuously on $f$.

The space $\hbox{Diff}^2(M)$ is a Banach manifold, see
\cite{abrahamrobbins} for details. For $f\in \hbox{Diff}^2(M)$,
the tangent space to $\hbox{Diff}^2(M)$ at $f$ has a natural
description. Let $X_f$ be the Banach space of $C^2$ sections of
the pullback bundle ${f}^* TM$, that is, the bundle whose fiber
over $p\in M$ is $T_{fp}M$. We write $ f + g $ to indicate the
diffeomorphism $\operatorname{exp}_{f} \circ g$. That is, if $g$
is a small vector field in $X_f$, then
$$(f + g)(p) = \operatorname{exp}_{f} (g(p))$$ is in $\hbox{Diff}^2(M)$ and
is close to $f$. So a small disk in $X_f$ is a chart for a small
neighborhood of $f$, and $X_f$ is thereby identified with the
tangent space $T_f\hbox{Diff}^2(M)$.

Define the map
\begin{equation*}
\begin{split}
\operatorname{Eval} : \mathcal{F}
\times M &\rightarrow \mathcal{F}
\times  M
\\
(f, p) &\mapsto (f, fp).
\end{split}
\end{equation*}
$\operatorname{Eval}$ is $C^2$ because left-composition is a
smooth operation on functions.

\begin{Lemm}\label{l.newcenter} If the diameter of ${\mathcal F}$ is
sufficiently
small, then $Eval$ has a partially hyperbolic splitting:
$$T({\mathcal F} \times M) = \mathbb{E}^u \oplus \mathbb{E}^c \oplus
\mathbb{E}^s,$$
where
$${\mathbb E}^u_{f,p} = 0\times E^u_p(f)\subset 0\times TM$$
$${\mathbb E}^s_{f,p} =  0 \times E^s_p(f) \subset 0\times TM,$$
and ${\mathbb E}^c_{f,p}$ is the graph of a linear map:
$$P_{f,p} : X_f \oplus  E^c_p(f)\to E^u_p(f)\oplus E^s_p(f).$$
\end{Lemm}

\begin{proof}
The tangent to $\operatorname{Eval}$ at $(f, p)$ acts on a vector
$\begin{bmatrix}
  g    \\
  v
\end{bmatrix} \in  T_{f,
p}(\mathcal{F} \times M)$ as
$$
T_{f,p}\operatorname{Eval}
\begin{bmatrix}
  g    \\
  v
\end{bmatrix}
=
\begin{bmatrix}
  g    \\
  g(p) + T_pf(v)
\end{bmatrix}
=
\begin{bmatrix}
  \operatorname{Id}_{\mathcal{F}}
&  0 \\
  \operatorname{ev}_{p}  &  T_pf
\end{bmatrix}
\begin{bmatrix}
  g    \\
  v
\end{bmatrix},
$$
where $\operatorname{ev}_{p}  $ evaluates the section of $f^{*}TM$
at $p$. In particular, this implies that the subbundles ${\mathbb
E}^u = 0 \times  E^u$, ${\mathbb E}^s = 0 \times  E^s$ cited above
are $T\operatorname{Eval}$-invariant. (The bundle $0 \times E^c$
is also $T\operatorname{Eval}$-invariant, but it is too small to
be the ${\mathbb E}^c$ we want.)

Note  that the subbundle $T{\mathcal F} \oplus 0$ is not
$T\operatorname{Eval}$-invariant, nor   is the subbundle
$T{\mathcal F} \oplus E^c$ whose fiber at $(f, p)$ is $X_f \oplus
E^c_p(f)$.  For if $v \in E^c_p(f)$ then the
$T\operatorname{Eval}$-image of $\begin{bmatrix}
  g \\
  v
\end{bmatrix}
$ is $
\begin{bmatrix}
  g  \\
  g(p) + T^cf_p(v)
\end{bmatrix}
$, and this vector need not lie in $T{\mathcal F} \oplus E^c$.
Nevertheless, by the domination hypotheses, the
$T\operatorname{Eval}$ graph transform defines a fiber contraction
of the bundle  whose fiber at $(f,p)$ is
$$
L(X_f \oplus (E^u\oplus E^c)_p(f), E^s(f)).
$$
The resulting invariant section is the unique
$T\operatorname{Eval}$-invariant subbundle ${\mathbb
E}^c\oplus{\mathbb E}^s \subset
  T(\mathcal{F}
\times  M)$ whose fiber at $(f,p)$ projects isomorphically onto
$X_f \oplus  (E^c\oplus E^s)_p(f)$.

Similarly, we find the unique
$T\operatorname{Eval}^{-1}$-invariant subbundle ${\mathbb
E}^c\oplus\mathbb{E}^u \subset
  T(\mathcal{F}
\times  M)$ whose fiber at $(f,p)$ projects isomorphically onto
$X_f \oplus  (E^c\oplus E^u)_p(f)$. Intersecting  these bundles, we
obtain the $T\operatorname{Eval}^{-1}$-invariant subbundle
${\mathbb E}^c$.
\end{proof}

\begin{Rmk} At the end of this section, we give a series expression for
$\mathbb{E}^{cu}$.
\end{Rmk}

\begin{Cor}\label{c.dDdE} Suppose that $f:M\to M$ is $C^2$ and
partially hyperbolic, with splitting:
$$TM = E^u\oplus E^c\oplus E^s.$$
If $f_0$ satisfies the pointwise bunching condition:
\begin{equation}
\begin{split}\label{e.bunchhyp}
  \sup_p
\frac{\norm{T^c_pf_0}} {\boldsymbol{m}(T^u_pf_0)
\boldsymbol{m}(T^c_pf_0)}  <  1,
\end{split}
\end{equation}
then, for all $p\in M$,  $E^u$ is continuously differentiable at
$(f_0,p)$ with respect to ${\mathbb E}^c$, where ${\mathbb E}^c$
is given by Lemma~\ref{l.newcenter}.
\end{Cor}

In fact this corollary can be stated in a more general form that
can be useful in applications.  Not only is it possible to
differentiate $E^u$ along ${\mathbb E}^c$, but in
fact bundles in dominated decompositions can be differentiated
along  ${\mathbb E}^c$ as well.  If $f_0$ has a dominated
decomposition
$$TM = R\oplus S,$$
then standard graph-transform arguments apply to show that
for $f$ sufficiently $C^1$-close to $f_0$, this decomposition
has a unique continuation
$$TM = R(f)\oplus S(f)$$
that is dominated for $Tf$.  Under appropriate bunching
hypotheses, we can differentiate $R(f)$ in the ${\mathbb E}^c$
direction:

\begin{Cor}\label{c.dDdE2} Suppose that $f:M\to M$ is $C^2$ and
partially hyperbolic, with splitting:
$$TM = E^u\oplus E^c\oplus E^s.$$
Suppose also that
$$TM = R\oplus S$$
is a dominated decomposition for $f_0$.
If $f_0$ satisfies the pointwise bunching condition:
\begin{equation}
\begin{split}\label{e.bunchyp2}
  \sup_p
\frac{\norm{T^S_pf_0}} {\boldsymbol{m}(T^R_pf_0)
\boldsymbol{m}(T^c_pf_0)}  <  1,
\end{split}
\end{equation}
then, for all $p\in M$,  $R$ is continuously differentiable at
$(f_0,p)$ with respect to ${\mathbb E}^c$, where ${\mathbb E}^c$
is given by Lemma~\ref{l.newcenter}.
\end{Cor}

\begin{proof}[Proof of Corollary~\ref{c.dDdE2}]
We first construct the bundle over ${\mathcal F}\times M$ whose
fiber over $(f,p)$ is the space of linear maps $L(R_p(f),
S_p(f))$.  Since $\operatorname{Eval}$ preserves the
factors $\{f\}\times M$, its tangent map $T\operatorname{Eval}$
induces a graph transform map on this bundle, covering
$\operatorname{Eval}$, which is a fiber contraction, with:
$$k_{f,p} = \frac{\norm{T^S_pf}} {\boldsymbol{m}(T^R_pf)}.$$
The unique invariant section of this graph transform is ${\mathbb
R} = 0\oplus R\subset T{\mathcal F}\times TM$. (note that
the bundle ${\mathbb R}$ is not to be confused with the real
numbers ${\bf R}$).

  Now suppose
$\gamma$ is any curve tangent to ${\mathbb E}^c$. As in the proof
of Theorem~\ref{t.sections}, we obtain differentiability of
${\mathbb R}$ (and hence, of $R$)  along $\gamma$ when
$k_{f,p}$ dominates the contraction along $\gamma$ at $(f,p)$  The
contraction along $\gamma$ at $(f,p)$ is bounded below by the
conorm of $T_{f,p}^{{\mathbb E}^c }\operatorname{Eval}$, which in
turn is approximately given by
$$\boldsymbol{m}(T_{f_0,p}^{{\mathbb E}^c }\operatorname{Eval}) =
\min \{1, \boldsymbol{m}(T^{c}
f_0)\}.$$

Hence, ${\mathbb R} = 0\oplus R$ is differentiable along
$\gamma$ if
\begin{equation}
\begin{split}\label{e.ratio}
  \sup_{f,p}
\frac{k_{f,p}}{ \min\{\boldsymbol{m}(T^c_pf_0),1\}} <  1.
\end{split}
\end{equation}
Since $k_{f,p} < 1$ for all p,
and by the bunching hypothesis (\ref{e.bunchyp2}),  $k_{f,p}
<\boldsymbol{m}(T^c_pf_0)$,
the condition in (\ref{e.ratio}) is satisfied.

Now we apply Theorem~\ref{t.sections} and conclude that
$\mathbb{R}$ is continuously differentiable along
$\mathbb{E}^c$.

Note that Theorem~\ref{t.sections} needs to be re-proved in this
more general context, but because its original proof relied on
uniform estimates (this was the only necessity for   the
compactness assumption on $M$), it is not hard to do.

\end{proof}

We are now ready to prove Theorem B.  As mentioned in the
introduction, Theorem B is a corollary of the following more general
result.

\begin{Thm}\label{t.functions}
Let $\{f_t:M\to M\}_{t\in(-\eps,\eps)}$ be a $C^2$ family of
$C^2$,  partially hyperbolic diffeomorphisms having,
for each $t\in(-\eps,\eps)$, a $Tf_t$-invariant splitting:
$$TM = E^u(f_t)\oplus E^c(f_t) \oplus E^s(f_t).$$

Then there exists $\epsilon_0>0$ so that, for every $p\in M$ and
every $v\in E^c(p)$,
there exists a $C^1$ path
$$\varphi_{p}:(-\epsilon_0,\epsilon_0)\to M$$
with the following properties:
\begin{enumerate}
\item $\varphi_{p,v}(0) = p$
\item $\dot{\varphi}_{p,v}(0)\in  v + E^u\oplus E^s$,
\item If
$$TM = R(f_0)\oplus S(f_0)$$
is any dominated decomposition for $f_0$
satisfying the pointwise bunching condition
\begin{equation}
\begin{split}
  \sup_p
\frac{\norm{T^S_pf_0}} {\boldsymbol{m}(T^R_pf_0)
\boldsymbol{m}(T^c_pf_0)}  <  1,
\end{split}
\end{equation}
then $t\mapsto R_{\varphi_{p,v}(t)} (f_t)$ is $C^1$.
\end{enumerate}
\end{Thm}

\begin{Rmk}
$E^c$ is allowed to be the trivial bundle in
Theorem~\ref{t.functions}, in which case $f_0$ is Anosov. If $f_0$
is Anosov, then ${\mathbb E}^c$ is uniquely integrable,
$\varphi_{p,0}$ is unique, and $p\mapsto \varphi_{p,0}(t)$ is the
homeomorphism conjugating $f_0$ to $f_t$.
\end{Rmk}

\begin{Rmk} Similarly, if $E^c$ is integrable and tangent to a
plaque-expansive  foliation ${\mathcal W}^c$, then ${\mathbb E}^c$
is also tangent to a foliation ${\mathbb W}^c$.
  The maps $p\mapsto
\varphi_{p,0}(t)$ can be shown to be leaf conjugacies between
$f_0$ and $f_t$.

If $f_0$ is $r$-normally-hyperbolic:
\begin{equation}
\begin{split}
\|T^cf\|^r < \boldsymbol{m}(T^uf) \qquad & \|T^sf\| <
\boldsymbol{m}(T^cf)^r,
\end{split}
\end{equation}
then the leaves of ${\mathbb W}^c$ are $C^r$. In this case,
$t\mapsto \varphi_{p,v}(t)$ can also be chosen to be $C^{r}$.

If, in addition, the stronger center bunching condition:
\begin{equation}
\begin{split}
  \sup_p
\frac{\norm{T^c_pf_0}} {\boldsymbol{m}(T^u_pf_0)
\boldsymbol{m}(T^c_pf_0)^r}  <  1
\end{split}
\end{equation}
holds, then the $C^r$ Section Theorem implies that ${\mathbb E}^u$
is $C^r$ along the leaves of ${\mathbb W}^c$ (and so $t\mapsto
E^u_{\varphi_{p,v}(t)}(f_t)$ is also $C^r$).
\end{Rmk}

\begin{Rmk} A simple refinement of the proof shows that
both  $t\mapsto R_{\varphi_{p,v}(t)} (f_t)$ and
$t\mapsto \varphi_{p,v}(t)$ are $C^{1+\alpha}$, where
there is a bound on the $\alpha$-H{\"o}lder norm of
the $t$-derivative that is uniform in $p,v$.  The exponent
$\alpha$  is determined by several bunching conditions.

\end{Rmk}

\begin{proof}[Proof of Theorem~\ref{t.functions}]
Let  $\mathcal F$ and, for $(f,q)\in {\mathcal F}\times M$,
the linear map $P_{f,q} : X_f \oplus  E^c_q(f)\to E^u_q(f)\oplus
E^s_q(f)$ be given by Lemma~\ref{l.newcenter}, so that
$${\mathbb E}^c_{f,q} = \hbox{graph}(P_{f,q}).$$
Choose $\epsilon_0>0$ so that $f_t\in {\mathcal F}$, for all
$t\in(-\epsilon_0,\epsilon_0)$.

We identify the submanifold
$$\{f_t\times M\,\mid\, t\in (-\epsilon_0,\epsilon_0)\subset
{\mathcal F}\times M$$
with $(-\epsilon_0,\epsilon_0)\times M$ in the obvious way. In the
tangent bundle ${\bf R}\times TM$ to this manifold, $P_{f_t,q}$
and ${\mathbb E}^c_{f_t,q}$ have their counterparts $P_{t,q}:{\bf
R}\oplus E^c_{f_t}(q) \to (E^u\oplus E^s)_{f_t}(q)$ and ${\mathbb
E}^c_{t,q}$, where
$${\mathbb E}^c_{t,q} = \hbox{graph}(P_{t,q}).$$

Let $p$ and $v$ be given, and let $V$ be any continuous vector
field on $M$ with the properties:
$V(p) = v$ and $V(q)\in E^c(q)$, for all
$q\in M$. Define  a vector field $\Omega$ on
$(-\epsilon_0,\epsilon_0)\times M$ by:
$$\Omega(t,q) = \frac{\partial}{\partial t} + V(q) +
P_{t,q}(\frac{\partial}{\partial t} + V(q)).$$
Notice that $\Omega(t,q)\in {\mathbb E}^c_{t,q}$, for all
$(t,q)\in (-\epsilon_0,\epsilon_0)\times M$.  It follows from
the bunching hypothesis and Corollary~\ref{c.dDdE2} that $R$ is
differentiable along the
integral curves of $\Omega$.

Let $\hat{\varphi}_{p,v}$ be any integral curve of $\Omega$
through $(p,0)$.  Now $\varphi_{p,v}$ is defined to be the $M$
coordinate of $\hat\varphi_{p,v}$: $$\hat{\varphi}_{p,v} = (t,
\varphi_{p,v}).$$ It is straightforward to check that
$\varphi_{p,v}$ satisfies (1)-(3).
\end{proof}

\subsection{A series expansion for
$\mathbb{E}^{cu}$}

  We give a series expression for
$\mathbb{E}^{cu}$ as follows. Define  the linear map $P^{cu}_{f,p}
: X \oplus E^{cu}_p(f) \rightarrow E^{s}_p(f)$ as the series
$$
P^{cu}_{f,p} (g, v) =
\sum_{k=0}^{\infty}
T^s_{f^{-k}p}f^k(g^s(f^{-k}p)).
$$
Note that the series does not depend
on $v$.   The domination conditions
imply that the series converges.  Under
$T\operatorname{Eval}$, the graph
of $P^{cu}_{f,p}$ is sent to the
graph of $P^{cu}_{f,fp}$.  Hence,
by uniqueness,
$$
\mathbb{E}^{cu}_{f,p} =
\operatorname{graph}
(P^{cu}_{f,p}).
$$
In the same way we get a unique
$T\operatorname{Eval}$-invariant
subbundle $\mathbb{E}^{cs} \subset
T(\mathcal{F} \times  M)$ whose
fiber at $(f, p)$  projects
isomorphically onto
$X \oplus E^{cs}_p(f)$, and
$
\mathbb{E}^{cs}_{f,p} =
\operatorname{graph}
(P^{cs}_{f,p})
$
where
$$
P^{cs}_{f,p}(g,v)
=
\sum_{k=0}^{\infty}
T^u_{f^kp}f^{-k}(g^u(f^kp)).
$$

The intersection of these two
subbundles is the center bundle
$\mathbb{E}^c$.  Namely, at $(f,
p)$, the fiber of $\mathbb{E}^c$ is
the graph of the map $P^c_{f, p} :
X \oplus E^c_p(f) \rightarrow
(E^u_p(f) \oplus E^s_p(f))$, where
$$
P^c_{f,p} (g,v)
=
\sum_{k=0}^{\infty}
T^u_{f^kp}f^{-k}(g^u(f^kp))
+
\sum_{k=0}^{\infty}
T^s_{f^{-k}p}f^k(g^s(f^{-k}p)).
$$

%%%%%%%%%%%%%%%%%%%%%%%%%%%%%%%%%%%%%%%%%%%%%

\section{When $f\mapsto E^u_p(f)$ actually {\em is}
differentiable}\label{s.Euisdiff}

We described in the previous section  how $f\mapsto E^u_p(f)$ is generally
not
differentiable, even if $f$ is Anosov.  In fact, if  $p\mapsto
E^u_p(f_0)$ fails to be differentiable in even one direction
at $p_0$,  then  $f\mapsto E^u_{p_0}(f)$ is
not differentiable at $f_0$.  For in that case, it is easy to construct a
smooth $1$-parameter family of diffeomorphisms $\varphi_t:M\to M$
such that $t\mapsto E^u_{\varphi_tp_0}(f_0)$ is not differentiable
at $t=0$; but then $E^u_{p_0}(\varphi_t f_0 \varphi_t^{-1})
= T\varphi_t(E^u_{\varphi_t p_0}(f_0))$
is not differentiable at $t=0$.

It turns out that, under the usual center bunching hypothesis,
nonsmoothness of $p\mapsto E^u_p(f_0)$ is the {\em only}
obstruction to differentiability of $f\mapsto E^u_p(f)$ at $f_0$
in all directions.

The results that follow apply to $1$-parameter families of $C^2$
diffeomorphisms
$\{f_t\}_{t\in I}$ such that $t\mapsto f_t$ is a $C^1$ map from $I$ into
$\hbox{Diff}^2(M)$ -- a $C^1$ family of $C^2$ diffeomorphisms.
Since the original proof of Theorem C is somewhat lengthy and the result
is subsumed by Theorem D, we omit the proof of Theorem C and
present instead a proof of Theorem D, following closely the approach of
Dolgopyat in \cite{D}.

Assume that for each $t  \in I$,
$f_t$ is partially hyperbolic with splitting
$$
TM = E_t^u \oplus E_t^c \oplus E_t^s.
$$
Write
$$
E_t = E_t^c  \quad H_t = E^u_t \oplus
E^s_t.
$$

%%--- THEOREM ----
\begin{Thm}[Theorem D]
\label{t: }If $E_0$ is a
$C^1$ bundle then the curve of
subbundles
$t
\mapsto  E_t$ is  differentiable
with respect to $t$ at $t =  0$,
and the derivative    $(
dE_t(x)/dt )_{t=0}$
   depends
continuously on
$x
\in  M$.
\end{Thm}
%%  -- EndTheorem

  %%%REMARK

\begin{Rmk}
Theorem D  remains valid, and
the proof is the same,  if the partially
hyperbolic splitting is replaced by a
dominated triple splitting
$R_t \oplus S_t
\oplus T_t$.  Namely,  the middle bundle
$S_t$  is differentiable
  with respect to
$t$ at
$t = 0$, provided  that $S_{x,0}$
is $C^1$.   Similarly, there is
nothing special about the
one-dimensionality of the
parameter $t$.
\end{Rmk}

The following facts about weak
continuity  will be used. We
assume that $W$ is a Banach
space, but that $W$ also carries a
weak topology.  Of course, if
$W$ has finite dimension,  the
weak and strong topologies
coincide. We have in
mind
  the case that
$W$ is a space of operators  on
the the infinite dimensional
Banach space of continuous
sections of a vector bundle and
$ \Lambda  =
\mathbb{R}$.

\begin{Defn}
A function $h : \Lambda
\rightarrow  W$ is \textbf{weakly
continuous}  at
$\mu \in \Lambda $
   if
$h(\lambda )$ tends weakly to
$h(\mu )$ and
$\norm{h(\lambda )}$ stays
bounded as $\lambda
\rightarrow  \mu $.
\end{Defn}

  %%%%--  Proposition
\begin{Prop}[Weak Inversion
Rule]
\label{p:weakinverse}
If a curve of invertible operators
$t \mapsto A_t$
   is   weakly continuous   at $t =
0$  and  if  the operators' conorms
are uniformly positive then the
curve of inverse operators is also
weakly continuous at $t = 0$.
\end{Prop}
%%  --- EndProposition

%%---  Proof --
\begin{proof}
Let $t \mapsto  A_t$ be the
curve of operators, and let $V$
be the Banach space on which
they operate.  Then, as $t
  \rightarrow  0$, $A_t$
converges weakly to $A_0$  and
$\norm{A_t - A_0}$ stays
bounded.  The conorm
assumption means that  for all
small
$t$, $\norm{A_t^{-1}} \leq  M$.

For each $v \in  V$,
\begin{equation*}
\begin{split}
\abs{A_t^{-1}(v) - A_0^{-1}(v)
}&=
\abs{A_t^{-1}\circ (A_0 - A_t) \circ
A_0^{-1}(v)}
\\
&\leq  M\abs{v -
A_t(A_0^{-1}(v))}.
\end{split}
\end{equation*}
  Since
$A_0^{-1}(v)$ is fixed, and $A_t$
converges weakly to $A_0$,
$A_t(A_0^{-1}(v)) \rightarrow v$ as $t
\rightarrow 0$, which completes the proof
that $ A_t^{-1}$ converges
weakly to $A_0^{-1}$ as
$t \rightarrow  0$.  But also,
\begin{equation*}
\begin{split}
\abs{A_t^{-1}(v) - A_0^{-1}(v)
}&=
\abs{A_t^{-1}\circ (A_0 - A_t) \circ
A_0^{-1}(v)}
\\
&\leq M\norm{A_0 - A_t}M\abs{v}
\end{split}
\end{equation*}
implies that $\norm{A_t^{-1} -
A_0^{-1}}$   stays bounded as $t
\rightarrow 0$, and completes
the proof that the inverse curve
is weakly continuous.
\end{proof}
%%  --- QED ---

\bigskip

Now we return to the  splitting
$TM = E_t
\oplus H_t$, where $H_t$ is the
hyperbolic part of the partially
hyperbolic splitting for
$f_t$, and $E_t$ is the  center
part.  We are assuming that
$E = E_0$ is a $C^1$ bundle.

Let $\widetilde{H}$ be a smooth
approximation to $H_0$, and
express
$Tf_t$ with respect to the
splitting $TM = E \oplus
\widetilde{H}$ as
$$
T_xf_t =
\begin{bmatrix}
  A_{x,t}  &  B_{x,t} \\
  C_{x,t}  &  K_{x,t}
\end{bmatrix}.
$$
Since $f_t$ is a $C^1$ curve of
$C^{2}$ diffeomorphisms,
$A, B, C, K$ are
$C^1$ functions of $x, t$.
At
$t = 0$ we have
$$
C_{x,0} = 0 \qquad \textrm{and} \qquad
A_{x,0} = T_xf_0|_{E}
$$
for all $x$.  Furthermore, when
$\widetilde{H}$ closely approximates
$H$,  $\norm{B}$ is small.
Consequently, if $P :  E \rightarrow
\widetilde{H}$ has norm $\leq  1$ then $A +
BP$ is invertible and the norm of its inverse
is uniformly bounded.  Uniformity refers to
$P, x, t$.

Let $\mathcal{L}$ be the vector
bundle over $M$ whose fiber at
$x$ is $\mathcal{L}_x = L(E_x,
\widetilde{H}_x)$.
Equipping
$\mathcal{L}_x$ with the
operator norm gives
$\mathcal{L}$  a Finsler; let
$\mathcal{L}(1)$ be its
unit ball bundle.  Denote by
$\operatorname{Sec}(\mathcal{
L})$ the Banach space of
continuous sections $X : M
\rightarrow
\mathcal{L}$, equipped with the
sup norm $\norm{\;\;}$.  Its unit
ball is
$\operatorname{Sec}
(\mathcal{L}(1))$.

$Tf_t$ defines a graph transform
$$
\begin{CD}
$$
  \mathcal{L}(1)
@>\text{\normalsize
$\qquad (Tf_t)_{\#} \qquad$}>>
\mathcal{L}\\
  @V\text{\normalsize
$\pi $}VV
@VV\text{\normalsize
$\pi $}V\\
M @>\text{\normalsize $\qquad
f_t \qquad$}>> M
$$
\end{CD}
$$
according to the condition
$T_xf_t
(\operatorname{graph}P) =
\operatorname{graph}(
(T_xf_t)_{\#}(P))$.
That is,
$$
(T_xf_t)_{\#}(P)
=
(C_{x,t} + K_{x,t}P)\circ
(A_{x,t} + B_{x,t}P)^{-1},
$$
which is a linear map $E_{f_tx}
\rightarrow
\widetilde{H}_{f_tx}$.  The
graph transform
naturally induces a nonlinear map
on the space of sections, %
\begin{equation*}
\begin{split}
G_t : \operatorname{Sec}(
\mathcal{L}(1)) &\rightarrow
\operatorname{Sec}(\mathcal{L})
\end{split}
\end{equation*}
such that
$$
G_t (X)  =  (Tf_t)_{\#} \circ X
\circ f_t^{-1}.
$$

%%%%--  Proposition
\begin{Prop}
\label{p:Gtanalytic}
$G_t$ is uniformly analytic.
\end{Prop}
%%  --- EndProposition

%%%REMARK

\begin{Rmk}
$(Tf_t)_{\#}$ is not analytic, it is
only $C^1$.  Nevertheless, for
each fixed $t$, its action on the
space of continuous sections is
analytic.   The uniformity refers to
$t$.
\end{Rmk}

  We prove
Proposition~\ref{p:Gtanalytic}
by factoring
$G_t$ into a product of
several analytic maps.
Let $\mathcal{E}$,
$\mathcal{E}_t$ and
$\mathcal{E}_t^{-1}$ denote the
bundles  whose fibers at $x \in
M$ are $\mathcal{E}_x =
L(E_x, E_x)$,
$\mathcal{E}_{x, t} = L(E_x,
E_{f_tx})$, and
$\mathcal{E}_{x,t}^{-1} =
L(E_{f_tx}, E_{x})$.
  Let $\mathcal{A}$,
$\mathcal{A}_t$, and
$\mathcal{A}_t^{-1}$ denote
the invertible elements in
$\mathcal{E}$,
$\mathcal{E}_t$ and
$\mathcal{E}_t^{-1}$,   and
denote
  sectional inversion as
$\operatorname{Inv} :
\operatorname{Sec}
(\mathcal{A})
\rightarrow
\operatorname{Sec}
(\mathcal{A})
$, $\operatorname{Inv}_t :
\operatorname{Sec}
(\mathcal{A}_t)
\rightarrow
\operatorname{Sec}(
\mathcal{A}_t^{-1})$.

%%-  Lemma -----
\begin{Lemm}
\label{l:sinversion}
Sectional inversion is uniformly
analytic.
\end{Lemm}
%%  -- EndLemma --

%%---  Proof --
\begin{proof}
Consider the identity section
$\operatorname{Id}$ of
$\mathcal{A}$.  Any section
near $\operatorname{Id}$ is
inverted by the power series
$$
A^{-1} = \sum_{k=0}^{\infty}
(\operatorname{Id} - A)^k,
$$
and hence sectional inversion is
analytic in a neighborhood of the
identity section. For $A$ in a
neighborhood of the  general
section
$A_0 : M \rightarrow
\mathcal{A}$, sectional
inversion factors according to
the commutative diagram
$$
\begin{CD}
$$
  A @>\text{\normalsize
$\operatorname{Inv}$ near
$A_0$ }>> A^{-1}\\
  @V\text{\normalsize
$L_{A_0^{-1}}$}VV
@AA\text{\normalsize
$R_{A_0}$}A\\
  A_0^{-1}A
@>\text{\normalsize
$\operatorname{Inv}$ near
$Id$}>> A^{-1}A_0
$$
\end{CD}
$$
where $L_{A_0^{-1}}$ and
$R_{A_0}$ are left and right
multiplication by the sections
$A_0^{-1}$ and $A_0$.  Since
$L_{A_0^{-1}}$ and
$R_{A_0}$   are continuous linear
transformations of the section
spaces, they are analytic, which
completes the proof of the lemma
for sections in a neighborhood of
the identity section.   The
corresponding diagram
$$
\begin{CD}
$$
  \operatorname{Sec}(
\mathcal{A}_t)
@>\text{\normalsize
$\operatorname{Inv}_t$ near
$A_0$ }>>
\operatorname{Sec}(
\mathcal{A}_t^{-1})\\
  @V\text{\normalsize
$L_{A_0^{-1}}$}VV
@AA\text{\normalsize
$R_{A_0}$}A\\
  \operatorname{Sec}(
\mathcal{A})
@>\text{\normalsize
$\operatorname{Inv}$ near
$Id$}>>
\operatorname{Sec}(
\mathcal{A})
$$
\end{CD}
$$
applies to sectional
inversion in the neighborhood of
a section
$A_0 : M \rightarrow
\mathcal{A}_t$, and shows that
$\operatorname{Inv}_t$ is
analytic.

Uniform analyticity  means that
for any
$r$, the $r^{\textrm{th}}$
   derivative of
$\operatorname{Inv}_t$  is
uniformly bounded on sets of
sections such that
$\norm{A}$ and
$\norm{A^{-1}}$ are uniformly
bounded; this is clear from the
higher order chain rule and the
factorization of sectional
inversion given above.
\end{proof}
%%  --- QED ---

%%---  Proof --
\begin{proof}
[Proof of
Proposition~\ref{p:Gtanalytic}]
We have $G_t(X) = (Tf_t)_{\#}
\circ X \circ f_t^{-1}$, and must
show that $G_t$ is a
uniformly analytic  function
of
$X \in
\operatorname{Sec}
(\mathcal{L})$.  We factor
$G_t$ as the Cartesian product
of two affine maps on section
spaces, followed by inversion in
one of the two spaces, followed
by sectional linear composition,
all of which is expressed by
commutativity of
$$
\begin{CD}
$$
\operatorname{Sec}
(\mathcal{L})
@>\text{\normalsize $\quad
\qquad
G_t \quad
\qquad $}>>
\operatorname{Sec}
(\mathcal{L})
\\
  @V\text{\normalsize
$\operatorname{Aff}_1
\times
\operatorname{Aff}_2$}VV
@AA\text{\normalsize
composition}A\\
  \operatorname{Sec}
(\mathcal{L}_t)
\times
\operatorname{Sec}
(\mathcal{A}_t)
@>\text{\normalsize $\qquad
\operatorname{Id}
\times
\operatorname{Inv}_t
  \qquad$}>>
\operatorname{Sec}(
\mathcal{L}_t)
\times
\operatorname{Sec}
(\mathcal{A}_t^{-1})
$$
\end{CD}
$$
where $\mathcal{L}_t$ is the
bundle over $M$ whose fiber at
$x$ is $L(E_x,
\widetilde{H}_{f_tx})$, and
$$
\operatorname{Aff}_1(X)
=
C_t + K_tX
\qquad
\operatorname{Aff}_2(X)
= A_t + B_tX.
$$
Uniform analyticity of $G_t$
then follows from
Lemma~\ref{l:sinversion}.
\end{proof}
%%  --- QED ---

  The
$r^{\textrm{th}}$-order
Taylor expansion of $G_t$ at the
zero section is
$$
G_t(X) =
Z_t + Q_t(X) +  \dots
+ \frac{1}{r!}(D^rG_t)_0(X^r) +
R_t(X) ,
$$
  where $Z_t = G_t(0)$, $Q_t =
(DG_t)_0$.

%%%%--  Proposition
\begin{Prop}
\label{p:Gtr}
For small $t$,
\begin{itemize}
   \item[(a)] $t \mapsto Z_t$ is
$C^1$.
   \item[(b)] $t \mapsto  (I -
Q_t)^{-1}$ is
weakly continuous.
   \item[(c)] $\norm{R_t(X)}/
\norm{X}^{2} $ is
uniformly bounded for all
small $X
\in
\operatorname{Sec}(
\mathcal{L})$.
\end{itemize}
\end{Prop}
%%  --- EndProposition

%%---  Proof --
\begin{proof}
At the zero
section, the
$0^{\textrm{th}}$ and first
derivatives of
$$
G_t(X) = (C_t + K_tX)(A_t +
B_tX)^{-1}\circ  f_t^{-1},
$$
with respect to $X$ are
computed at once as
\begin{equation*}
\begin{split}
Z_t &= (C_tA_t^{-1})\circ
f_t^{-1}
\\
Q_t(X) &= (K_tXA_t^{-1} +
C_tA_t^{-1}B_tXA_t^{-1})
\circ  f_t
^{-1}
\end{split}
\end{equation*}

Since $f_t$ is a $C^1$ curve of
$C^{2}$ diffeomorphisms, and
since the splitting $E \oplus
\widetilde{H}$ is $C^1$, the
curves $t \mapsto A_t$, $t
\mapsto B_t$, $t \mapsto  C_t$,
$t \mapsto K_t$ in the
appropriate bundles are
$C^1$.  This makes (a)
immediate, and also shows that
the curve $t \mapsto  Q_t$ in
$\operatorname{Sec}(\mathcal{L})$
is weakly continuous.

By inspection, at $t = 0$,
$Q_t$ becomes the
hyperbolic operator
$$
Q_0(X) = (K_0XA_0^{-1})\circ
f_0^{-1},
$$
because $C_{t=0} = 0$.  Thus,
for all small $t$,    $I
- Q_t$ is uniformly invertible,
and
Proposition~\ref{p:weakinverse}
implies that $t \mapsto
(I-Q_t)^{-1}$ is weakly
continuous.

Assertion (c) follows from the
Mean Value Theorem and the
fact that the
second derivative
of $G_t$ is uniformly bounded
near the zero section.
\end{proof}
%%  --- QED ---

%%---  Proof --
\begin{proof}[Proof of Theorem
D]
Proposition~\ref{p:Gtr} implies
that
$$
G_t(X) = Z_t + Q_t(X)
+ R_t(X)
$$
and $\norm{R_t(X)} =
O(1)\norm{X}^2$ as $\norm{X}
\rightarrow  0$. Let
$P_t : x \mapsto P_{x,t}$ be the
unique
$G_t$-invariant section of
$ \mathcal{L}  $ with norm $\leq
  1$.  Thus $P_{x,t} : E_x
\rightarrow
\widetilde{H}_x$ and
$$
E_{x,t} = \operatorname{graph}P_{x, t}
=
\{v + P_{x,t}(v) \in  T_xM :  v \in
E_x\}.
$$
Theorem D asserts that $E_t$ is
differentiable
at $t = 0$.   That is,
$$
\frac{dP_{x,t}}{dt}\Big|_{t=0}
$$
exists and   is
continuous with respect to $x$.

Plugging   $X = P_t $ into
the  Taylor expansion of $G_t$
gives
$$
P_t = Z_t + Q_t(P_t) +
R_t(P_t),
$$
and since $I - Q_t$ is invertible, we get
$$
P_t =  (I - Q_t)^{-1}(Z_t +
R_t(P_t)).
$$
Thus
\begin{equation}
\begin{split}
\label{e:Ptquadratic}
  \norm{P_t}  \leq  \norm{(I -
Q_t)^{-1}} (\norm{Z_t}  +
M\norm{P_t} ^2).
\end{split}
\end{equation}
(These norms refer to section
sup-norms or to operator norms,
as appropriate.)

   Now we estimate
$Z_t = (C_t
\circ  A_t^{-1}) \circ  f_t^{-1}$
as follows. It is differentiable with
respect to $t$, and since $C_{t=0}
= 0$, we have
$Z_{t=0} = 0$.  Thus
$\norm{Z_t} = O(1)t$
as $t \rightarrow 0$.  Since $P_t$
is continuous in $t$, and $P_0 =
0$, we get
$\norm{P_t}_0^2 \ll \norm{P_t}_0$ when
$t$ is small, which lets us absorb
the squared term into the l.h.s. of
the inequality
(\ref{e:Ptquadratic}), so
$$
\norm{P_t} = O(1)t
$$
as $t \rightarrow  0$.
Consequently, we get a bootstrap
effect on the remainder:
$$
\norm{R_t(P_t)} = O(1)t^2
$$
as $t \rightarrow 0$.
Combined with the
more exact estimate on
$Z_t$,
$$
  Z_t = t
Z^{\prime}_0 +
o(1)t
$$
where $Z_0^{\prime} =
(d/dt)_{t=0}(Z_t)$, this gives
$$
\frac{P_t}{t} =  (I-Q_t)^{-1}
Z^{\prime}_0 + (I-Q_t)^{-1}(o(1)
+ O(1)t).
$$
Proposition~\ref{p:Gtr} implies
that $(I - Q_t)^{-1}$ converges
weakly to $(I - Q_0)^{-1}$ as $t
\rightarrow  0$, so
$$
\lim_{t\rightarrow 0}
(I-Q_t)^{-1} Z^{\prime}_0
=
(I-Q_0)^{-1}  Z^{\prime}_0,
$$
while uniform boundedness of
$\norm{(I-Q_t)^{-1}}$  implies that
$$
\lim_{t\rightarrow 0} (I-Q_t)^{-1}(o(1)) +
O(1)t) = 0.
$$
  Thus, as $t \rightarrow 0$,
$$
\frac{P_{x,t} - P_{x,0}}{t}
\rightarrow  (I-Q_0)^{-1}
Z^{\prime}_0,
$$
uniformly in $x \in  M$, which
completes the proof that
$t \mapsto E_t$ is differentiable
at
$t = 0$, and that its derivative
there, $(I -
Q_0)^{-1}Z_0^{\prime}$,
depends continuously on
$x
\in  M$.
\end{proof}
%%  --- QED ---

  %%%REMARK

\begin{Rmk}
Suppose that $E_0$ and $Df_t$
are
$C^r$, $r \geq 2$.  We tried to
show that  $E_t$ is
$r^{\textrm{th}}$-order
differentiable at $t = 0$ in
the sense that there is an
$r^{\textrm{th}}$ order  Taylor
expansion for
$E_t$  at
$t = 0$.  Many ingredients of
the preceding proof of the $r =
1$ case above generalize very
nicely to $r  \geq 2$.  There is a
natural notion of weak
$r^{\textrm{th}}$-order
differentiability, and it behaves
well with respect to operator
inversion and operator
products.  However, we would
also need affirmative answers to
the following two questions:
\begin{itemize}
   \item[(a)] Is the curve $t
\mapsto  (I - Q_t)^{-1}$
in
$\operatorname{Sec}
(\mathcal{L})$ weakly
differentiable at $t = 0$?
   \item[(b)] Does the operator
$(I - Q_0)^{-1}$ send $C^1$
sections of $\mathcal{L}$ to
$C^1$ sections?
\end{itemize} At first, it would
be acceptable to assume
analyticity of
$E_0$ and $f_t$.

\end{Rmk}


\begin{thebibliography}{1}
\bibitem{abrahamrobbins} Abraham, Ralph, and Robbins, Joel
\emph{Transversal Mappings and Flows},
Benjamin, New York, 1967.
\bibitem{D} Dolgopyat, Dmitry, {On differentiability of SRB states,}
preprint.
\bibitem{GPS}  Grayson, Matthew; Pugh, Charles; Shub, Michael,
{\em Stably ergodic diffeomorphisms.}  Ann. of
Math. (2)  {\bf 140 } (1994), no. 2, 295--329.
\bibitem{HPS}  Hirsch, M. W.; Pugh, C. C.; Shub, M.
{\em Invariant manifolds.} Lecture Notes in Mathematics, Vol. 583.
Springer-Verlag, Berlin-New
York, 1977.
\bibitem{KKPW} Katok, A.; Knieper, G.; Pollicott, M.; Weiss, H.
{\em Differentiability and analyticity of topological entropy for
Anosov and geodesic
flows.} Invent. Math. {\bf 98} (1989), no. 3, 581--597.
\bibitem{LMM} de la Llave, R.; Marco, J. M.; Moriyon, R.
{\em Canonical perturbation theory of Anosov
systems and regularity results for the Liv\v sic cohomology
equation.} Ann. of Math. (2) {\bf 123 }
(1986), no.
3, 537--611.
\bibitem{M} Ma{\~n}{\'e}, Ricardo,
{\em On the dimension of the compact invariant sets of certain nonlinear
maps.}
Dynamical systems and turbulence, Warwick 1980 (Coventry, 1979/1980),
pp. 230--242,
Lecture Notes in Math., {\bf 898}, Springer, Berlin-New York, 1981.
\bibitem{Milnor}  Milnor, John, {\em Fubini foiled: Katok's
paradoxical example in measure theory. }
Math. Intelligencer  {\bf 19 }(1997), no. 2, 30--32.
\bibitem{PS} Pugh, Charles; Shub, Michael,
{\em Stably ergodic dynamical systems and partial hyperbolicity.} J.
Complexity {\bf 13} (1997), no. 1, 125--179.
\bibitem{PS2}  Pugh, Charles; Shub, Michael,
{\em  Stable ergodicity and julienne quasi-conformality.}  J. Eur.
Math. Soc. (JEMS) {\bf 2} (2000), no. 1, 1--52.
\bibitem{ruelle1}  Ruelle, David,
{\em Differentiation of SRB states.}  Comm. Math. Phys. {\bf 187
}(1997), no. 1, 227--241.
\bibitem{ruelle2} Ruelle, David, ``Extensions of a result by Shub and
Wilkinson..." preprint
\bibitem{SW} Shub, Michael; Wilkinson, Amie,
{\em Pathological foliations and removable zero exponents.} Invent.
Math. {\bf 139} (2000), no. 3, 495--508.
\bibitem{W} Wilkinson, Amie {\em Stable ergodicity of the time-one
map of a geodesic flow.}
Ergodic
Theory Dynam. Systems {\bf 18} (1998), no. 6, 1545--1587.
\end{thebibliography}
\end{document}